# A Transformation-Proximal Bundle Algorithm for Multistage Adaptive Robust Optimization and Application to Constrained Robust Optimal Control

Chao Ning, Fengqi You*

Robert Frederick Smith School of Chemical and Biomolecular Engineering, Cornell University, Ithaca, New York 14853, USA
*Corresponding author. Phone: (607) 255-1162; Fax: (607) 255-9166; E-mail: fengqi.you@cornell.edu

**Abstract**

This paper presents a novel transformation-proximal bundle algorithm for multistage adaptive robust optimization problems. By partitioning recourse decisions into state and control decisions, the proposed algorithm applies affine control policy only to state decisions and allows control decisions to be fully adaptive, thus transforming the original problem into an equivalent two-stage Adaptive Robust Optimization (ARO) problem. Importantly, this multi-to-two transformation is general enough to be employed with any two-stage ARO solution algorithms, thus opening a new avenue for a variety of multistage ARO algorithms. The proximal bundle method is developed for the resulting two-stage problem along with convergence analysis. In an inventory control application, the affine disturbance-feedback control policy suffers from a severe suboptimality with an average gap of 34.88%, while the proposed algorithm generates an average gap of merely 1.68%.

## 1. Introduction

In recent years, robust optimization has become an increasingly popular methodology to immunize optimization problems against uncertainty among both control and optimization communities (Petersen & Tempo, 2014; Wu et al., 2013). Robust optimization can be roughly classified into three categories: static robust optimization, two-stage Adaptive Robust Optimization (ARO), and multistage ARO. Two-stage ARO allows recourse decisions to be adaptive to realized uncertainties (Ben-Tal et al., 2004), thus typically generating less conservative solutions than static robust optimization. To overcome the limitation of two-stage structures, multistage ARO emerges as a practical yet more computationally challenging paradigm for sequential decision making, such as finite-horizon optimal control.

This paper proposes a novel multi-to-two transformation-proximal bundle algorithm. In a multistage decision-making setting, decision variables can be partitioned into two different groups, namely state decision variables and control/local decision variables (Bodur & Luedtke, 2017; Ning & You, 2018a). We propose a novel multi-to-two transformation scheme that converts the multistage ARO problem into an equivalent two-stage counterpart. Specifically, by enforcing only state decisions to be affine functions of uncertainty, Multistage Adaptive Robust Mixed-Integer Linear Program (MARMILP) is reduced into a Two-stage Adaptive Robust Mixed-Integer Linear Program (TARMILP). We prove that such transformation is valid if state decisions follow causal control policies. The multi-to-two transformation scheme is general enough to be combined with existing two-stage ARO solution algorithms for solving MARMILPs.

The remainder of this paper is organized as follows. Section 2 presents a novel multi-to-two transformation scheme for multistage ARO problems. A transformation-proximal bundle algorithm is then proposed in Section 3. We further develop a lower bounding technique in Section 4. In Section 5, an application to robust optimal inventory control is presented, followed by a concluding remark.

## 2. The multi-to-two transformation scheme

The MARMILP in its general form is shown as follows:

$$\min_{\mathbf{x},\, \mathbf{s}_t(\cdot),\, \mathbf{y}_t(\cdot)} \max_{\mathbf{u}\in U} \mathbf{c}'\mathbf{x} + \sum_{t=1}^{T}\left[\mathbf{d}_t'\mathbf{s}_t(\mathbf{u}^t) + \mathbf{f}_t'\mathbf{y}_t(\mathbf{u}^t)\right]$$

$$\text{s.t. } \mathbf{T}_t\mathbf{x} + \mathbf{A}_t\mathbf{s}_t(\mathbf{u}^t) + \mathbf{B}_t\mathbf{s}_{t-1}(\mathbf{u}^{t-1}) \quad (1)$$
$$+ \mathbf{W}_t\mathbf{y}_t(\mathbf{u}^t) = \mathbf{h}_t^0 + \mathbf{H}_t\mathbf{u}_t, \forall \mathbf{u}\in U, t$$
$$\mathbf{L}_t\mathbf{x} + \mathbf{E}_t\mathbf{s}_t(\mathbf{u}^t) + \mathbf{G}_t\mathbf{y}_t(\mathbf{u}^t) \geq \mathbf{m}_t^0 + \mathbf{M}_t\mathbf{u}_t, \ \forall \mathbf{u}\in U, t$$

The material in this paper was partially, in particular the material of section 2 and some of the material of section 3, presented at the 57th IEEE Conference on Decision and Control, December 17–19, 2018, Miami Beach, FL, USA. The authors acknowledge financial support from the National Science Foundation (NSF) CAREER Award (CBET-1643244).

where $T$ is the total number of time stages, $\mathbf{u}_1, \ldots, \mathbf{u}_T$ are uncertainties revealed over $T$ stages, $\mathbf{x}$ is a vector of "here-and-now" decisions made prior to any uncertainty realizations, $\mathbf{s}_1, \ldots, \mathbf{s}_T$ are adjustable state decision variables, and $\mathbf{y}_1, \ldots, \mathbf{y}_T$ are adjustable control decision variables. Note that the "here-and-now" decisions and adjustable decisions involve continuous decision variables. The prime symbol ′ stands for the transpose of a generic vector. Let vector $\mathbf{u}^t=[\mathbf{u}_1', \ldots, \mathbf{u}_t']'$ be the concatenated vectors of past uncertainty realizations, and $\mathbf{u}=[\mathbf{u}_1', \ldots, \mathbf{u}_T']'$. $\mathbf{c}$, $\mathbf{d}_t$ and $\mathbf{f}_t$ are the vectors of cost coefficients.

The proposed multi-to-two transformation only restricts state $\mathbf{s}_t(\cdot)$ to follow an affine control policy as shown in (2).

$$\mathbf{s}_t(\mathbf{u}^t) = \mathbf{P}_t \mathbf{u}^t + \mathbf{q}_t \quad (2)$$

where $\mathbf{P}_t$ and $\mathbf{q}_t$ are the coefficients of the affine function. After plugging the control policy (2) into (1), the MARMILP can be formulated as follows:

$$\min_{\substack{\mathbf{x}, \mathbf{P}_t, \mathbf{q}_t \\ \mathbf{y}_t(\cdot)}} \max_{\mathbf{u} \in U} \left( \mathbf{c}'\mathbf{x} + \sum_{t=1}^{T} \mathbf{d}_t' \mathbf{q}_t \right) + \sum_{t=1}^{T} \left[ \mathbf{d}_t' \mathbf{P}_t \mathbf{u}^t + \mathbf{f}_t' \mathbf{y}_t(\mathbf{u}^t) \right]$$
$$\text{s.t.} \ \mathbf{A}_t \left( \mathbf{P}_t \mathbf{u}^t + \mathbf{q}_t \right) + \mathbf{B}_t \left( \mathbf{P}_{t-1} \mathbf{u}^{t-1} + \mathbf{q}_{t-1} \right) \quad (3)$$
$$+ \mathbf{W}_t \mathbf{y}_t(\mathbf{u}^t) = \mathbf{h}_t^0 + \mathbf{H}_t \mathbf{u}_t - \mathbf{T}_t \mathbf{x}, \ \forall \mathbf{u} \in U, t$$
$$\mathbf{E}_t \left( \mathbf{P}_t \mathbf{u}^t + \mathbf{q}_t \right) + \mathbf{G}_t \mathbf{y}_t(\mathbf{u}^t) \geq \mathbf{m}_t^0 + \mathbf{M}_t \mathbf{u}_t - \mathbf{L}_t \mathbf{x}, \ \forall \mathbf{u} \in U, t$$

where control decision $\mathbf{y}_t(\cdot)$ is a general function of uncertainty realizations.

For the ease of exposition, we present the nested formulation of multistage ARO problem (3) in (4).

$$\min_{\hat{\mathbf{x}} \in \Omega_0} \left\{ f_0(\hat{\mathbf{x}}) + \max_{\mathbf{u}_1 \in U_1} \min_{\mathbf{y}_1 \in \Omega_1(\hat{\mathbf{x}}, \mathbf{u}^1)} \left\{ f_1(\mathbf{y}_1) + \cdots \max_{\mathbf{u}_T \in U_T} \min_{\mathbf{y}_T \in \Omega_T(\hat{\mathbf{x}}, \mathbf{u}^T)} f_T(\mathbf{y}_T) \right\} \right\} \quad (4)$$

where $\hat{\mathbf{x}} = \{\mathbf{x}, \mathbf{P}_t, \mathbf{q}_t\}$ is an aggregated "here-and-now" decisions, set $\Omega_0$ represents its feasible region, and set $\Omega_t(\hat{\mathbf{x}}, \mathbf{u}^t)$ is the feasible region of adjustable control decisions at stage $t$ as given in (5). $U_1, \ldots, U_T$ denote uncertainty sets.

$$\Omega_t(\hat{\mathbf{x}}, \mathbf{u}^t) = \left\{ \mathbf{y}_t \left| \begin{array}{c} \mathbf{A}_t(\mathbf{P}_t \mathbf{u}^t + \mathbf{q}_t) + \mathbf{B}_t(\mathbf{P}_{t-1}\mathbf{u}^{t-1} + \mathbf{q}_{t-1}) \\ + \mathbf{W}_t \mathbf{y}_t = \mathbf{h}_t^0 + \mathbf{H}_t \mathbf{u}_t - \mathbf{T}_t \mathbf{x} \\ \mathbf{E}_t(\mathbf{P}_t \mathbf{u}^t + \mathbf{q}_t) + \mathbf{G}_t \mathbf{y}_t \geq \mathbf{m}_t^0 + \mathbf{M}_t \mathbf{u}_t - \mathbf{L}_t \mathbf{x} \end{array} \right. \right\} \quad (5)$$

The objective functions in the nested multistage ARO formulation (4) are explicitly defined in (6).

$$\begin{cases} f_0(\hat{\mathbf{x}}) = \mathbf{c}'\mathbf{x} + \sum_{t=1}^{T} \mathbf{d}_t' \mathbf{q}_t \\ f_t(\mathbf{y}_t) = \mathbf{f}_t' \mathbf{y}_t + \mathbf{d}_t' \mathbf{P}_t \mathbf{u}^t, \ t=1, \ldots, T \end{cases} \quad (6)$$

The uncertainty set for stage $t$ is given by (7),

$$U_t = \text{Proj}_{\mathbf{u}_t} \left( U \big|_{\mathbf{u}_1, \ldots, \mathbf{u}_{t-1}} \right) \quad (7)$$

where $U_t$ is defined as the projection of uncertainty set $U$ onto $\mathbf{u}_t$ given the values of $\mathbf{u}_1$ to $\mathbf{u}_{t-1}$.

**Theorem 1.** If the affine control policy (2) is used only for adjustable state decisions, the multistage ARO problem (1) is transformed into a two-stage ARO problem given below.

$$\min_{\hat{\mathbf{x}} \in \Omega_0} \left\{ \left( \mathbf{c}'\mathbf{x} + \sum_{t=1}^{T} \mathbf{d}_t' \mathbf{q}_t \right) + \max_{\mathbf{u} \in U} \min_{\mathbf{y} \in \{\mathbf{y} | \mathbf{y}_t \in \Omega_t(\hat{\mathbf{x}}, \mathbf{u}^t), \forall t\}} \sum_{t=1}^{T} \left( \mathbf{d}_t' \mathbf{P}_t \mathbf{u}^t + \mathbf{f}_t' \mathbf{y}_t \right) \right\} \quad (8)$$

where $\mathbf{y}=[\mathbf{y}_1', \ldots, \mathbf{y}_T']'$ be the concatenated control decisions.
**Proof.** Considering the max-min optimization problem in (4) at $t=T-1$, we have

$$\max_{\mathbf{u}_{T-1} \in U_{T-1}} \min_{\mathbf{y}_{T-1} \in \Omega_{T-1}(\hat{\mathbf{x}}, \mathbf{u}^{T-1})} \left\{ f_{T-1}(\mathbf{y}_{T-1}) + \max_{\mathbf{u}_T \in U_T} \min_{\mathbf{y}_T \in \Omega_T(\hat{\mathbf{x}}, \mathbf{u}^T)} f_T(\mathbf{y}_T) \right\}$$
$$= \max_{\mathbf{u}_{T-1} \in U_{T-1}} \max_{\mathbf{u}_T \in U_T} \left\{ \min_{\mathbf{y}_{T-1} \in \Omega_{T-1}(\hat{\mathbf{x}}, \mathbf{u}^{T-1})} f_{T-1}(\mathbf{y}_{T-1}) + \min_{\mathbf{y}_T \in \Omega_T(\hat{\mathbf{x}}, \mathbf{u}^T)} f_T(\mathbf{y}_T) \right\} \quad (9)$$
$$= \max_{\{\mathbf{u}_{T-1}, \mathbf{u}_T\} \in \text{Proj}_{\{\mathbf{u}_{T-1}, \mathbf{u}_T\}}\left(U|_{\mathbf{u}_1, \ldots \mathbf{u}_{T-2}}\right)} \left\{ \sum_{t=T-1}^{T} \min_{\mathbf{y}_t \in \Omega_t(\hat{\mathbf{x}}, \mathbf{u}^t)} f_t(\mathbf{y}_t) \right\}$$

The derivation can be performed backward until $t=1$.

$$\min_{\hat{\mathbf{x}} \in \Omega_0} \left\{ f_0(\hat{\mathbf{x}}) + \max_{\{\mathbf{u}_1, \ldots, \mathbf{u}_T\} \in \text{Proj}_{\{\mathbf{u}_1, \ldots, \mathbf{u}_T\}}(U)} \left\{ \sum_{t=1}^{T} \min_{\mathbf{y}_t \in \Omega_t(\hat{\mathbf{x}}, \mathbf{u}^t)} f_t(\mathbf{y}_t) \right\} \right\}$$
$$= \min_{\hat{\mathbf{x}} \in \Omega_0} \left\{ f_0(\hat{\mathbf{x}}) + \max_{\mathbf{u} \in U} \left\{ \sum_{t=1}^{T} \min_{\mathbf{y}_t \in \Omega_t(\hat{\mathbf{x}}, \mathbf{u}^t)} f_t(\mathbf{y}_t) \right\} \right\} \quad (10)$$
$$= \min_{\hat{\mathbf{x}} \in \Omega_0} \left[ f_0(\hat{\mathbf{x}}) + \max_{\mathbf{u} \in U} \min_{\mathbf{y} \in \{\mathbf{y} | \mathbf{y}_t \in \Omega_t(\hat{\mathbf{x}}, \mathbf{u}^t), \forall t\}} \sum_{t=1}^{T} f_t(\mathbf{y}_t) \right]$$

According to (6) and (10), the multistage problem (4) is equivalent to problem (8), which concludes the proof. □

## 3. Transformation-proximal bundle algorithm

### 3.1. A multistage robust optimization solution algorithm

The worst-case recourse function of the two-stage ARO problem, denoted as $Q(\hat{\mathbf{x}})$, is shown in (11).

$$Q(\hat{\mathbf{x}}) = \max_{\mathbf{u} \in U} \min_{\mathbf{y} \in \{\mathbf{y} | \mathbf{y}_t \in \Omega_t(\hat{\mathbf{x}}, \mathbf{u}^t), \forall t\}} \sum_{t=1}^{T} \left( \mathbf{d}_t' \mathbf{P}_t \mathbf{u}^t + \mathbf{f}_t' \mathbf{y}_t \right) \quad (11)$$

where the "max-min" optimization problem is often referred to as an adversarial optimization problem.

Problem (8) can be considered as a minimization problem whose objective function is given by (12).

$$F(\hat{\mathbf{x}}) = \left( \mathbf{c}'\mathbf{x} + \sum_{t=1}^{T} \mathbf{d}_t' \mathbf{q}_t \right) + Q(\hat{\mathbf{x}}) \quad (12)$$

where $F(\hat{\mathbf{x}})$ is the objective function of two-stage ARO problem (8).

Due to the multi-level optimization structure, the objective function $F(\hat{\mathbf{x}})$ does not have an analytical expression. In the proximal bundle method, bundle information includes the past query points $\hat{\mathbf{x}}^l$ ($l=1, \ldots, k$), their corresponding function values $F(\hat{\mathbf{x}}^l)$, and sub-gradients of function $F$ at these query points. We need to solve the max-min optimization problem in (11) to obtain the function value and a sub-gradient at one query point. To this end, the two-level optimization problem (11) is equivalently reformulated into sub-problem (**SUP**) using KKT condition.

$$\max_{\substack{\mathbf{u} \in U \\ \mathbf{y}_t, \boldsymbol{\varphi}_t, \boldsymbol{\pi}_t}} \sum_{t=1}^{T} \left( \mathbf{d}_t' \mathbf{P}_t \mathbf{u}^t + \mathbf{f}_t' \mathbf{y}_t \right)$$
$$\text{s.t. } \mathbf{W}_t' \boldsymbol{\varphi}_t + \mathbf{G}_t' \boldsymbol{\pi}_t = \mathbf{f}_t, \ \forall t \qquad \text{(SUP)}$$
$$\boldsymbol{\pi}_t \geq \mathbf{0}, \ \forall t$$
$$\mathbf{y}_t \in \Omega_t(\hat{\mathbf{x}}, \mathbf{u}^t), \ \forall t$$
$$\left[ \mathbf{G}_t \mathbf{y}_t + \mathbf{E}_t (\mathbf{P}_t \mathbf{u}^t + \mathbf{q}_t) - \mathbf{m}_t^0 - \mathbf{M}_t \mathbf{u}_t + \mathbf{L}_t \mathbf{x} \right]_i (\boldsymbol{\pi}_t)_i = 0, \ \forall t,i$$

where $\boldsymbol{\varphi}_t$ and $\boldsymbol{\pi}_t$ are the dual variables for the constraints in (5) at stage $t$, and $i$ denotes the element index of a vector.

We linearize these constraints into (13) by using the big-$M$ method, which is a standard technique (Wolsey, 1998).

$$(\boldsymbol{\pi}_t)_i \leq M \cdot (\mathbf{w}_t)_i, \ \forall t,i$$
$$\left[ \mathbf{G}_t \mathbf{y}_t + \mathbf{E}_t (\mathbf{P}_t \mathbf{u}^t + \mathbf{q}_t) - \mathbf{m}_t^0 - \mathbf{M}_t \mathbf{u}_t + \mathbf{L}_t \mathbf{x} \right]_i \leq M \cdot \left[ 1 - (\mathbf{w}_t)_i \right], \ \forall t,i \qquad (13)$$

where $\mathbf{w}_t$ represents a vector of binary decision variables. To ensure the exactness of big-$M$ reformulations, $M$ should be chosen as a sufficiently large positive number.

With sub-gradients and function values, we build the optimality cutting plane model for $F(\hat{\mathbf{x}}^l)$ shown in (14).

$$\tilde{F}_k(\hat{\mathbf{x}}) = \max_{l=1,\ldots,k} \left\{ F(\hat{\mathbf{x}}^l) + \langle \mathbf{g}^l, \hat{\mathbf{x}} - \hat{\mathbf{x}}^l \rangle \right\} \qquad (14)$$

where $\tilde{F}_k(\hat{\mathbf{x}})$ is the optimality cutting plane model at the $k$-th iteration. $\mathbf{g}^l$ is one sub-gradient of the objective function $F$ at the $l$-th query point and can be obtained using optimal dual variables in the same way as the Benders decomposition (Thiele et al., 2009). Note that the two-stage ARO problem (8) may not satisfy the relative complete recourse assumption. Therefore, for some query point $\hat{\mathbf{x}}^l$, there exist certain uncertainty realizations that render the second-stage optimization problem infeasible. This implies $F(\hat{\mathbf{x}}^l) = +\infty$ or equivalently $\hat{\mathbf{x}}^l \notin \text{dom } F$, where dom represents the domain of a function. For a given $\hat{\mathbf{x}}^l$, we either derive a lower linearization of function $F$ (optimality cut) or obtain a cutting plane that separates $\hat{\mathbf{x}}^l$ and dom $F$ (feasibility cut). To check whether $\hat{\mathbf{x}}^l \in \text{dom } F$ or not, the following Feasibility Problem (**FP**) needs to be solved.

$$\max_{\mathbf{u} \in U} \min_{\mathbf{y}_t, \boldsymbol{\alpha}_t^+, \boldsymbol{\alpha}_t^-, \boldsymbol{\beta}_t} \sum_t \left( \mathbf{1}' \boldsymbol{\alpha}_t^+ + \mathbf{1}' \boldsymbol{\alpha}_t^- + \mathbf{1}' \boldsymbol{\beta}_t \right)$$
$$\text{s.t. } \mathbf{A}_t (\mathbf{P}_t \mathbf{u}^t + \mathbf{q}_t) + \mathbf{B}_t (\mathbf{P}_{t-1} \mathbf{u}^{t-1} + \mathbf{q}_{t-1}) + \mathbf{W}_t \mathbf{y}_t + \boldsymbol{\alpha}_t^+ - \boldsymbol{\alpha}_t^-$$
$$= \mathbf{h}_t^0 + \mathbf{H}_t \mathbf{u}_t - \mathbf{T}_t \mathbf{x}, \ \forall t \qquad \text{(FP)}$$
$$\mathbf{E}_t (\mathbf{P}_t \mathbf{u}^t + \mathbf{q}_t) + \mathbf{G}_t \mathbf{y}_t + \boldsymbol{\beta}_t \geq \mathbf{m}_t^0 + \mathbf{M}_t \mathbf{u}_t - \mathbf{L}_t \mathbf{x}, \ \forall t$$
$$\boldsymbol{\alpha}_t^+, \boldsymbol{\alpha}_t^-, \boldsymbol{\beta}_t \geq \mathbf{0}, \ \forall t$$

where $\boldsymbol{\alpha}_t^+$, $\boldsymbol{\alpha}_t^-$, and $\boldsymbol{\beta}_t$ are slack variables, and $\mathbf{1}$ is the vector of ones in an appropriate dimension. Let $\omega(\hat{\mathbf{x}}^l)$ denote the optimal value of problem (**FP**) associated with a query point $\hat{\mathbf{x}}^l$. If $\omega(\hat{\mathbf{x}}^l) = 0$, there exist feasible second-stage decisions for any uncertainty realizations in uncertainty set $U$. Thus, we have $\hat{\mathbf{x}}^l \in \text{dom } F$ and only need optimality cuts. If $\omega(\hat{\mathbf{x}}^l) > 0$, the worst-case uncertainty realization can lead to the nonexistence of feasible recourse decisions. As a result, the feasibility cut is required.

To determine the next query point, we consider the Moreau-Yosida regularization of $\tilde{F}_k(\hat{\mathbf{x}})$ given by (15),

$$G(\hat{\mathbf{x}}) = \tilde{F}_k(\hat{\mathbf{x}}) + \frac{1}{2t_k} \|\hat{\mathbf{x}} - \mathbf{z}^k\|^2 \qquad (15)$$

where $\mathbf{z}^k$ is the stability center for the $k$-th iteration and $t_k$ is the positive proximal parameter (Hiriart-Urruty & Lemaréchal, 2013). We iteratively refine the cutting plane models by adding new query points on the fly. The optimal solution of Master Problem (**MP**) provides next query point.

$$\min_{\hat{\mathbf{x}} \in \Omega_0, \eta} \eta + \frac{1}{2t_k} \|\hat{\mathbf{x}} - \mathbf{z}^k\|^2$$
$$\text{s.t. } \eta \geq F(\hat{\mathbf{x}}^l) + \langle \mathbf{g}^l, \hat{\mathbf{x}} - \hat{\mathbf{x}}^l \rangle, \ l \in L_o \qquad \text{(MP)}$$
$$0 \geq F(\hat{\mathbf{x}}^l) + \langle \mathbf{g}_f^l, \hat{\mathbf{x}} - \hat{\mathbf{x}}^l \rangle, \ l \in L_f$$

where $\eta$ is an auxiliary variable. $L_o$ and $L_f$ denote the index sets of optimality and feasibility cuts, respectively. Akin to the Benders decomposition, constraint $\eta \geq F(\hat{\mathbf{x}}^l) + \langle \mathbf{g}^l, \hat{\mathbf{x}} - \hat{\mathbf{x}}^l \rangle$ corresponds to an optimality cut, while $0 \geq F(\hat{\mathbf{x}}^l) + \langle \mathbf{g}_f^l, \hat{\mathbf{x}} - \hat{\mathbf{x}}^l \rangle$ is a feasibility cut. Besides the cuts derived in the dual space, optimality cuts in the primal space can be added as well (Zeng & Zhao, 2013). In the proximal bundle method, the expected decrease $\delta_k$ defined in equation (16) is used to determine whether to update the stability center. Also, the expected decrease is used to check the stopping criterion. The full details of the algorithm can be found in (Ning & You, 2018b).

$$\delta_k = F(\mathbf{z}^k) - \tilde{F}_k(\hat{\mathbf{x}}^{k+1}) - \frac{1}{2t_k} \|\hat{\mathbf{x}}^{k+1} - \mathbf{z}^k\|^2 \qquad (16)$$

where $\hat{\mathbf{x}}^{k+1}$ is an optimal solution to (**MP**).

*3.2. Convergence analysis*

The proofs of propositions and lemmas in this section are given in (Ning & You, 2018b).

**Proposition 1.** $F(\hat{\mathbf{x}})$ in (12) is a convex function in $\hat{\mathbf{x}}$.

To facilitate exposition, we define the linearization errors at the stability center $\mathbf{z}^k$ in (17).

$$e_l = F(\mathbf{z}^k) - \left[ F(\hat{\mathbf{x}}^l) + \langle \mathbf{g}^l, \mathbf{z}^k - \hat{\mathbf{x}}^l \rangle \right], \ \forall l \qquad (17)$$

With the definition of $e_l$, we can rewrite $\tilde{F}_k(\hat{\mathbf{x}})$ in (18).

$$\tilde{F}_k(\hat{\mathbf{x}}) = \max_{l=1,\ldots,k} \left\{ F(\mathbf{z}^k) - e_l - \langle \mathbf{g}^l, \mathbf{z}^k - \hat{\mathbf{x}}^l \rangle + \langle \mathbf{g}^l, \hat{\mathbf{x}} - \hat{\mathbf{x}}^l \rangle \right\}$$
$$= F(\mathbf{z}^k) + \max_{l=1,\ldots,k} \left\{ -e_l + \langle \mathbf{g}^l, \hat{\mathbf{x}} - \mathbf{z}^k \rangle \right\} \qquad (18)$$

To prove the convergence of the proposed algorithm, we first present lemmas and proposition as follows (Belloni, 2005; Hiriart-Urruty & Lemaréchal, 2013).

**Lemma 1.** Consider the following Regularized Optimization Problem (ROP):

$$\min_{\hat{\mathbf{x}}} \tilde{F}_k(\hat{\mathbf{x}}) + \frac{1}{2t_k} \|\hat{\mathbf{x}} - \mathbf{z}^k\|^2 \qquad \text{(ROP)}$$

Then, the dual problem is shown in (19).

$$\max_{\left\{\boldsymbol{\alpha}\in R_+^k\left|\sum_{l=1}^k\alpha_l=1\right.\right\}} F(\mathbf{z}^k)-\frac{t_k}{2}\left\|\sum_{l=1}^k\alpha_l\cdot\mathbf{g}^l\right\|^2-\sum_{l=1}^k\alpha_l\cdot e_l \quad (19)$$

**Lemma 2.** Suppose $\boldsymbol{\alpha}$ is an optimal solution to the optimization problem in (19). Then, we have

(i) $\hat{\mathbf{g}}^k \in \partial \tilde{F}_k(\hat{\mathbf{x}}^{k+1})$ ; (ii) $\tilde{F}_k(\hat{\mathbf{x}}^{k+1}) = F(\mathbf{z}^k)-t_k\|\hat{\mathbf{g}}^k\|^2 - \hat{e}_k$ ;

(iii) $\delta_k = \frac{t_k}{2}\|\hat{\mathbf{g}}^k\|^2 + \hat{e}_k$ ; (iv) $\hat{\mathbf{g}}^k \in \partial_{\hat{e}_k} F(\mathbf{z}^k)$ .

where $\hat{\mathbf{g}}^k = \sum_{l=1}^k \alpha_l \cdot \mathbf{g}^l$ and $\hat{e}_k = \sum_{l=1}^k \alpha_l e_l$.

**Definition 1** (*Serious Steps*). For the proposed algorithm, serious steps refer to those steps in which the stability center is changed.

**Lemma 3.** Suppose $F^*$ be the optimal value of $\min F(\hat{\mathbf{x}})$ and $F^* > -\infty$. Then, we have inequality (20).

$$\sum_{k\in L_s} \delta_k \leq \frac{F(\mathbf{z}^0)-F^*}{m} < \infty \quad (20)$$

where $L_s$ denotes the set of iteration having serious steps.

**Assumption 1.** For the infinite number of serious steps, i.e. $|L_s|=+\infty$, the sequence $\{F(\mathbf{z}^k)\}_{k\in L_s}$ is assumed to converge and $\lim_{k\in L_s} F(\mathbf{z}^k) = F_* > -\infty$.

**Lemma 4.** For an infinite number of serious steps, we have

(i) If $\sum_{k\in L_s} t_k = \infty$, then $\liminf_{k\to\infty}\|\hat{\mathbf{g}}^k\|=0$; (ii) If $0<t_k\leq c$ and $\arg\min_{\hat{\mathbf{x}}} F(\hat{\mathbf{x}}) \neq \varnothing$, then $\{\mathbf{z}^k\}_{k\in L_s}$ is bounded.

**Definition 2** (*Null Steps*). Null steps are those steps in which the stability center remains the same.

**Lemma 5.** If there is a finite number of serious steps, i.e. $|L_s|<+\infty$, let $k_0$ be the index of last serious step, $\{\hat{\mathbf{x}}^k\}_{k\geq k_0}$ be the sequence of null steps, and $\mathbf{z}^{k_0}$ be the stability center generated by the last serious step. Then, we have ($k>k_0$),

$$F(\mathbf{z}^{k_0}) - \delta_k + \frac{1}{2t_k}\|\hat{\mathbf{x}} - \hat{\mathbf{x}}^{k+1}\|^2$$
$$= \tilde{F}_k(\hat{\mathbf{x}}^{k+1}) + \langle\hat{\mathbf{g}}^k, \hat{\mathbf{x}}-\hat{\mathbf{x}}^{k+1}\rangle + \frac{1}{2t_k}\|\hat{\mathbf{x}}-\mathbf{z}^{k_0}\|^2 \quad (21)$$

**Lemma 6.** For the proposed algorithm, the following equality and inequality hold.

(i) $\tilde{F}_k(\hat{\mathbf{x}}^{k+1}) + \langle\hat{\mathbf{g}}^k, \hat{\mathbf{x}}^{k+2}-\hat{\mathbf{x}}^{k+1}\rangle = F(\mathbf{z}^k)+\langle\hat{\mathbf{g}}^k, \hat{\mathbf{x}}^{k+2}-\mathbf{z}^k\rangle - \hat{e}_k$

(ii) $\tilde{F}_k(\hat{\mathbf{x}}^{k+1}) + \langle\hat{\mathbf{g}}^k, \hat{\mathbf{x}}^{k+2}-\hat{\mathbf{x}}^{k+1}\rangle \leq \tilde{F}_{k+1}(\hat{\mathbf{x}}^{k+2})$

**Assumption 2.** Sequence $\{t_k\}$ is nonincreasing and bounded below by a positive number.

**Propositon 2.** If there is a finite number of serious steps, let $k_0$ be the index of last serious step, $\{\hat{\mathbf{x}}^k\}_{k\geq k_0}$ the sequence of null steps, and $\mathbf{z}^{k_0}$ is the stability center generated by the last serious step. Then $\delta_k \to 0$ as $k\to\infty$.

**Theorem 2.** (i) For $\delta_{tol}=0$, there is a cluster point of sequence $\{\mathbf{z}^k\}$, denoted as $\mathbf{z}^\infty$, satisfying $F(\mathbf{z}^\infty)=F(\hat{\mathbf{x}}^*)$; (ii) for $\delta_{tol}>0$, the algorithm converges in finite steps.

**Proof.** For $\delta_{tol}=0$, the transformation-proximal bundle algorithm loops forever. There are two exclusive scenarios: (1) The algorithm implements an infinite number of serious steps; (2) After a finite number of serious steps, the algorithm implements only null steps. For scenario 1, we have $\hat{\mathbf{g}}^k \in \partial F_{\hat{e}_k}(\mathbf{z}^k)$ based on **Lemma 2**. Based on definition of approximate subgradient and the Cauchy-Schwarz inequality, we have $-\|\hat{\mathbf{g}}^k\|\cdot\|\hat{\mathbf{x}}^*-\mathbf{z}^k\| \leq F(\hat{\mathbf{x}}^*)-F(\mathbf{z}^k)+\hat{e}_k$. Since $\hat{\mathbf{x}}^*$ is an optimal solution, we further have (22).

$$F(\hat{\mathbf{x}}^*) \leq F(\mathbf{z}^k) \leq F(\hat{\mathbf{x}}^*) + \hat{e}_k + \|\hat{\mathbf{g}}^k\|\cdot\|\hat{\mathbf{x}}^*-\mathbf{z}^k\| \quad (22)$$

Note that $\|\hat{\mathbf{x}}^*-\mathbf{z}^k\|$ is bounded, because of **Lemma 3** and **Lemma 4**, $\lim_{k\to\infty}\hat{e}_k = 0$ and $\liminf_{k\to\infty}\|\hat{\mathbf{g}}^k\|=0$. Based on (22), there is a cluster point of sequence $\{\mathbf{z}^k\}$, denoted as $\mathbf{z}^\infty$, satisfying $F(\mathbf{z}^\infty)=F(\hat{\mathbf{x}}^*)$. Under scenario 2, we have $\delta_k \to 0$ as $k\to\infty$ based on **Proposition 2**. Also, $\delta_k \to 0$ implies $\hat{e}_k \to 0$ and $\|\hat{\mathbf{g}}^k\|\to 0$ as $k\to\infty$. according to **Lemma 2** (iii) and **Assumption 2**. Since stability center at $k$-th iteration remains at $\mathbf{z}^{k_0}$ for $k>k_0$, we have $F(\mathbf{z}^{k_0})=F(\hat{\mathbf{x}}^*)$ by taking the limit on (22). Thus, statement (i) still holds. For $\delta_{tol}>0$, suppose the algorithm does not converge in finite number of iterations, then we have $\delta_k > \delta_{tol} > 0, \forall k$, which contradicts $\delta_k \to 0$. based on **Lemma 3** and **Proposition 2**. □

## 4. The lower bounding technique

The idea of the proposed scenario-tree-based lower bounding approach is to replace the uncertainty set in MARMILPs with a finite number of uncertainty scenarios constructed within the cutting-plane algorithm. The Scenario-Tree Multistage Adaptive Robust Mixed-Integer Linear Program (STMARMILP), as shown in (23).

$$\begin{aligned}
\min_{\mathbf{x},\mathbf{s}_t,\mathbf{y}_t,\theta}\quad & \mathbf{c}'\mathbf{x}+\theta \\
\text{s.t.}\quad & \theta \geq \sum_{t=1}^T\left[\mathbf{d}_t'\mathbf{s}_t(\mathbf{u}_{(i)}^t)+\mathbf{f}_t'\mathbf{y}_t(\mathbf{u}_{(i)}^t)\right],\ \forall i\in\{1,\cdots,N\} \\
& \mathbf{A}_t\mathbf{s}_t(\mathbf{u}_{(i)}^t)+\mathbf{B}_t\mathbf{s}_{t-1}(\mathbf{u}_{(i)}^{t-1})+\mathbf{W}_t\mathbf{y}_t(\mathbf{u}_{(i)}^t) \\
& = \mathbf{h}_t^0+\mathbf{H}_t\mathbf{u}_{(i)}^t - \mathbf{T}_t\mathbf{x},\ \forall i\in\{1,\cdots,N\},\ t \\
& \mathbf{E}_t\mathbf{s}_t(\mathbf{u}_{(i)}^t)+\mathbf{G}_t\mathbf{y}_t(\mathbf{u}_{(i)}^t) \geq \mathbf{m}_t^0+\mathbf{M}_t\mathbf{u}_{(i)}^t - \mathbf{L}_t\mathbf{x},\ \forall i\in\{1,\cdots,N\},\ t \\
& \mathbf{u}_{(i)}^t = \mathbf{u}_{(j)}^t \Rightarrow \mathbf{s}_t(\mathbf{u}_{(i)}^t)=\mathbf{s}_t(\mathbf{u}_{(j)}^t),\ \mathbf{y}_t(\mathbf{u}_{(i)}^t)=\mathbf{y}_t(\mathbf{u}_{(j)}^t)\ \forall i,j,t
\end{aligned} \quad (23)$$

The optimality gap is defined $(UB-LB)/0.5(UB+LB)$, where $UB$ denotes the upper bound, and $LB$ represents the lower bound obtained via the STMARMILP.

**Theorem 3.** For any specific problem instance of MARMILPs, the following inequalities (24) hold.

$$v^S \leq v^* \leq v^{TPB} \leq v^{ADR} \quad (24)$$

where $v^S$, $v^*$, $v^{TPB}$, and $v^{ADR}$ present the optimal values of STMARMILP, MARMILP, TARMILP, and the affinely adjustable robust counterpart.

**Proof.** Since the scenario set is a subset of the uncertainty set ($\tilde{U} \subseteq U$), the scenario-tree counterpart STMARMILP is a relaxation of the original multistage ARO problem by satisfying only a subset of constraints. Hence, the objective value of STMARMILP provides a lower bound for the original multistage ARO problem ($v^S \leq v^*$).

In the original MARMILP, the recourse decisions are general functions of uncertainty. In both the affine control policy and the proposed transformation proximal bundle algorithm, all or some of the recourse variables are restricted to a fixed functional form of uncertainty realizations, thus providing upper bounds to the optimal value of the original multistage ARO problems ($v^* \leq v^{ADR}$ and $v^* \leq v^{TPB}$). Additionally, any feasible solution of the affinely adjustable robust counterpart is also feasible for the TARMILP. Therefore, we have $v^{ADR} \leq v^{TPB}$. □

## 5. Application: Robust optimal inventory control

In this section, we apply the proposed algorithm to robust finite-horizon optimal inventory control problems. We refer the readers to (Ning & You, 2018b) for more details and computational results of this application. Comparisons between the affine control policy (Ben-Tal, et al., 2004), the piecewise affine control policy (Chen & Zhang, 2009), and the proposed algorithm are made. In the application, the affine control policy suffers from severe suboptimality. Its largest relative gap can reach as high as 53.43%, and the average relative gap is 25.72%. By contrast, the control policy determined by the proposed algorithm has a relative gap of 1.33% on average, while its highest relative gap is merely 4.27%. We further compare the proposed lower bounding technique with a data-driven approach that samples uncertainty scenarios from the uncertainty set following the uniform distribution (Maggioni et al., 2016). To guarantee a fair comparison, the same number of uncertainty scenarios is used. As can be observed from **Fig. 1**, the proposed approach generates tighter lower bounds in each instance.

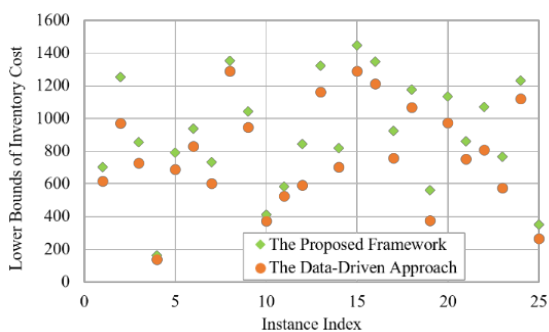

**Fig. 1.** Lower bounds of multi-period inventory cost determined by the proposed method and the data-driven approach.

To investigate the performance under different number of time stages, we implement computational experiments with $T$=10 and $T$=15. The average relative gap of affine control policy soars significantly from 25.72% to 34.88% when the value of $T$ changes from 5 to 15. In stark contrast, the average gap of the proposed control policy is increased by only 0.35%.

## 6. Conclusions

In this paper, a transformation-proximal bundle algorithmic framework was proposed. By employing the proposed multi-to-two transformation scheme, the multistage ARO problem was proved to be transformed into an equivalent two-stage ARO problem. The proximal bundle algorithm was further developed with the convergence analysis. The computational results demonstrated the effectiveness of the proposed control policy in robust inventory control.